\documentclass[letterpaper,conference]{ieeeconf}
\IEEEoverridecommandlockouts                              
\usepackage[utf8]{inputenc}
\usepackage{textcomp}

\usepackage{graphicx}
\usepackage{amsmath}
\usepackage[version=4]{mhchem}
\usepackage{siunitx}
\usepackage{graphicx}
\usepackage{nicefrac}
\usepackage{circuitikz}
\usepackage{float}
\usepackage{longtable,tabularx}
\usepackage[style=ieee ,sorting=none,maxbibnames=99,giveninits, doi=false]{biblatex}
\usepackage{tikz}
\usetikzlibrary{positioning}
\AtBeginBibliography{\small}
\usepackage{amsmath} 
\usepackage{amssymb}  
\usepackage{amsfonts}
\usepackage{amsthm}
\usepackage{pifont}
\newtheoremstyle{indented}
    {3pt}
    {3pt}
    {\addtolength{\leftskip}{0em}}
    {}
    {\bfseries}
    {.}
    {1em}
    {}
\theoremstyle{indented}

\newtheorem{proposition}{Proposition}[section]

\newtheorem{definition}{Definition}[section]
\newtheorem{fact}{Fact}[section]

\newtheorem{remark}{Remark}[section]

\usepackage{float} 
\usepackage[skip=1pt,font=footnotesize]{caption}

\usepackage{subcaption}

\usepackage{etaremune}

\usepackage{tikz}
\usetikzlibrary{shapes,arrows,calc,positioning}
\tikzstyle{bigblock} = [draw, fill=blue!20, rectangle, 
    minimum height=6em, minimum width=8em]
\tikzstyle{medblock} = [draw, fill=blue!20, rectangle, 
    minimum height=4em, minimum width=4em]    
\tikzstyle{mux} = [draw, fill=black!20, rectangle, 
    minimum height=5em, minimum width=0.1em]    
\tikzstyle{smallblock} = [draw, fill=blue!20, rectangle, 
    minimum height=2em, minimum width=3em]
    
\tikzstyle{data_block} = [draw, fill=green!20, rectangle, 
    minimum height=2em, minimum width=3em]
\tikzstyle{ops_block} = [draw, fill=blue!20, rectangle, 
    minimum height=2em, minimum width=3em]    
\tikzstyle{est_block} = [draw, fill=red!20, rectangle, 
    minimum height=2em, minimum width=3em]    
    
\tikzstyle{sum} = [draw, fill=blue!20, circle, node distance=1cm,minimum height=0.5cm]
\tikzstyle{signal} = [coordinate]
\tikzstyle{pinstyle} = [pin edge={to-,thin,black}]
\tikzstyle{block} = [draw, fill=blue!20, rectangle, 
    minimum height=3em, minimum width=9em]
\tikzstyle{blockS} = [draw, fill=blue!20, rectangle, 
    minimum height=3em, minimum width=4em]    
\tikzstyle{input} = [coordinate]
\tikzstyle{output} = [coordinate]
\usetikzlibrary{matrix}

\usetikzlibrary{positioning}
\usetikzlibrary{math}

\usepackage{hyperref}
\usepackage{xcolor}
\hypersetup{
    colorlinks,
    linkcolor={blue!100!black},
    citecolor={blue!50!black},
    urlcolor={blue!80!black}
}

\newcommand{\bc}{\begin{center}}
\newcommand{\ec}{\end{center}}
\newcommand{\benum}{\begin{enumerate}}
\newcommand{\eenum}{\end{enumerate}}
\newcommand{\nn}{\nonumber}
\newcommand{\matl}{\left[ \begin{array}}
\newcommand{\matr}{\end{array} \right]}
\newcommand{\matls}{\left[ \begin{smallmatrix}}
\newcommand{\matrs}{\end{smallmatrix} \right]}
\newcommand{\isdef}{\stackrel{\triangle}{=}}

\newcommand{\inv}{^{-1}}

\newcommand{\tr}{{\rm tr}\,}

\newcommand{\rmH}{{\rm H}}

\newcommand{\rmT}{{\rm T}}

\newcommand{\rmd}{{\rm d}}

\newcommand{\rmf}{{\rm f}}

\newcommand{\rmi}{{\rm i}}

\newcommand{\rmp}{{\rm p}}

\newcommand{\BBC}{{\mathbb C}}

\newcommand{\BBR}{{\mathbb R}}

\renewcommand{\matl}{\begin{bmatrix}}
\renewcommand{\matr}{\end{bmatrix} }

\usepackage{comment}
\usepackage{marginnote}

\title{Adaptive Quantum Control}

\title{Model-free, Learning-based Quantum Control System}
\title{Model-free Adaptive Control of LGKS Quantum System}
\title{Adaptive Fidelity-Based Density Tracking for Open Quantum Systems 
}

\author{
Jhon Manuel Portella Delgado
and
Ankit Goel
\thanks{Jhon Manuel Portella Delgado is a graduate student in the Department of Mechanical Engineering, University of Maryland, Baltimore County, 1000 Hilltop Circle, Baltimore, MD 21250. {\tt\small jportella@umbc.edu}}%
\thanks{Ankit Goel is an Assistant Professor in the Department of Mechanical Engineering, University of Maryland, Baltimore County,1000 Hilltop Circle, Baltimore, MD 21250. {\tt\small ankgoel@umbc.edu }}%
}

\begin{document}

\maketitle

\begin{abstract}
This paper presents an online learning-based adaptive control framework for density-matrix tracking in a two-level Lindblad--Gorini--Kossakowski--Sudarshan (LGKS) quantum system, in which the feedback control law does not require prior knowledge of the system Hamiltonian or dissipative operators. 
The adaptive controller is based on a continuous-time formulation of retrospective cost adaptive control (RCAC). 
To preserve the geometric structure of the quantum-state evolution, an adaptive PID controller driven by \textit{Uhlmann’s fidelity} is employed. 
The proposed approach is validated in numerical simulations for both low-entropy and high-entropy density-tracking tasks, and robustness to measurement noise in the feedback path is investigated.


\end{abstract}
\textit{\bf keywords:} 
Quantum feedback control, open quantum systems, online learning, adaptive PID control, Lindblad dynamics

\section{INTRODUCTION}

Quantum dynamical systems have several applications in sensing, metrology, communication, computation, and the internet \cite{giovannetti2011advances, toth2014quantum, taylor2016quantum, cacciapuoti2020entanglement, nielsen2001quantum}. 
%
The reliable operation of quantum circuits, such as quantum logic gates, requires the precise and predictable behavior of quantum dynamical systems, much like how classical digital circuits depend on the accurate design of their fundamental electrical components, such as resistors and transistors \cite{barenco1995elementary}.
In quantum computing, the basic unit of information is the qubit, which, unlike classical bits limited to two distinct states ($0$ or $1$), can exist in a superposition of both states, enabling quantum algorithms to perform computations exponentially faster than classical systems \cite{grover1996fast,shor1999polynomial}.
The successful realization of quantum devices thus depends on the precise control of quantum dynamical systems.
%

Quantum control concerns the manipulation of quantum systems to achieve specific dynamical or state-transfer objectives and plays a central role in quantum computing, communication, and sensing. 
Such control is typically realized through externally applied fields, for example, electromagnetic control inputs, that steer the system’s evolution. 
The desired quantum state depends on the application. 
In quantum computing, metrology, communication, and optical systems \cite{yamamoto2024engineering, zaletel2021preparation}, low-entropy (high-purity) states are often desirable to preserve coherence and reduce uncertainty, whereas in other settings, high-entropy or mixed states may be of interest when intrinsic randomness or thermal effects are central to the application. 
Examples of the latter include heat engines based on spin systems, quantum cryptography, and quantum thermodynamics simulations \cite{von2019spin, Sehgal2021, somhorst2023quantum}. 
Quantum control strategies typically seek to shape the system’s time evolution, governed by quantum operators such as Hamiltonians and density matrices. 
Closed quantum systems evolve according to the Schr\"odinger equation, whereas open quantum systems interacting with an environment are commonly modeled using the Lindblad--Gorini--Kossakowski--Sudarshan (LGKS) equation \cite{brif2010control}. 
Recent years have seen significant progress in control techniques for both isolated and open quantum systems \cite{dong2010quantum, cong2014control}.

In practice, the performance of model-based quantum control strategies is limited by uncertainty in the underlying system dynamics, unmodeled environmental effects, and constraints on the accuracy and availability of quantum state measurements. 
These challenges motivate the development of adaptive quantum control methods that adjust controller parameters online based on system feedback to compensate for modeling errors, parameter drift, and disturbances \cite{Kosut2003, kosut2013adaptive, Liu2020}. 
By incorporating online learning mechanisms, adaptive control can mitigate the impact of incomplete system knowledge and imperfect measurements \cite{egger2014adaptive}. 
Such approaches are relevant for practical quantum technologies, including quantum error correction and quantum information processing, where uncertainties and environmental interactions directly affect performance.

Despite these advances, most existing adaptive quantum control approaches primarily address parametric uncertainty while still requiring partial knowledge of the system model, such as the Hamiltonian structure or dissipative operators, which may be difficult to obtain accurately for complex quantum systems. 
Consequently, much of the existing literature focuses on robustness to parameter variations rather than control design in the absence of an explicit dynamical model \cite{egger2014adaptive, Koswara2014, Hu2016}. 
This motivates the investigation of control architectures in which the feedback law itself does not rely on detailed knowledge of the system dynamics.

This paper considers the problem of designing an adaptive controller for an open quantum system in which the control law does not require prior knowledge of the system Hamiltonian or dissipative operators. 
Specifically, we develop a continuous-time adaptive controller based on the Retrospective Cost Adaptive Control (RCAC) framework \cite{ali2015retrospective, Goel2015, Sobolic2016, poudel2023learning} for density-matrix tracking in open quantum systems governed by the Lindblad--Gorini--Kossakowski--Sudarshan (LGKS) equation \cite{manzano2020short}. 
In this preliminary study, Uhlmann’s fidelity is adopted as the scalar performance metric driving the adaptation \cite{UHLMANN1976273}, providing a geometrically meaningful measure of similarity between quantum states. 
The main contributions of this paper are:
\begin{itemize}
    \item the formulation of a continuous-time RCAC-based adaptive controller for LGKS quantum dynamics without requiring a plant model in the control law,
    \item the design of an adaptive PID feedback structure driven by fidelity-based error, and
    \item numerical validation of the proposed approach for tracking both low-entropy and high-entropy target states in a two-level open quantum system.
\end{itemize}

The remainder of the paper is organized as follows. 
Section \ref{sec:problem_formulation} describes the LGKS equation in the context of a controlled open quantum system. 
Section \ref{sec:CTRCAC} presents the mathematical foundations of continuous-time RCAC with forgetting factors. 
Section \ref{sec:numerical_simulation} presents numerical simulations validating the proposed adaptive controller on a two-level quantum system. 
Finally, Section \ref{sec:conclusions} concludes the paper and outlines directions for future work.

\section{Problem Formulation}
\label{sec:problem_formulation}

We consider the Lindblad-Gorini-Kossakowski-Sudarshan system as described in \cite{manzano2020short}, which is governed by
\begin{align}
    {\dot{\rho}}
        &=
            -i\hbar[H,
            {{\rho}}
            ] + \sum_{i=1}^m\left(L_i
            {{\rho}}
            L_i^\rmH  - \dfrac{1}{2}\{L_i^\rmH L_i, 
            {{\rho}}
            \}\right),
        \label{eq:Lindblad_equation}
\end{align}
\normalsize
where 
${{\rho}} \in \mathbb{C}^{n \times n}$
is the density matrix, 
the system Hamiltonian $H \in \mathbb{C}^{n \times n}$ is Hermitian, 
$m$ is the number of damping terms, and 
for $i \in (1, 2, \ldots, m), $ the matrices $L_i \in \BBC^{n \times n}$ are the jump operators.
Without loss of generality, the units are chosen such that $\hbar = 1.$
The diagonal elements of 
${{\rho}}$ 
represent the probability of a quantum event. 
It is shown in \cite[p~123]{breuer2002theory} that,
if 
$\tr {{\rho}}(0) = 1,$ 
then, for all $t > 0,$ 
${{\rho}}(t)$ 
given by \eqref{eq:Lindblad_equation} satisfies 
$\tr ({{\rho}}(t)) = 1.$
Furthermore, for all $t \geq 0,$ 
${{\rho}}(t)$ 
is positive semidefinite and Hermitian. 
%



In this paper, we consider a \textit{two-level quantum dynamic system}, that is, $n=2$ and thus ${{\rho}} \in \BBC^{2 \times 2}.$
As shown in \cite{altafini2012modeling}, the Hamilotonian $H$ can be decomposed as 
\begin{align}
    H
        &=
            H_0 + H_1u, 
        \label{eq:H_controlled}
\end{align}
where $H_0 \in \mathbb{C}^{2 \times 2}$ is the free Hamiltonian, 
$H_1 \in \mathbb{C}^{2 \times 2}$ is the control Hamiltonian, and 
$u \in \mathbb{R}$ is the scalar control input. 
The restriction $u \in \mathbb{R}$ ensures that $H$ remains Hermitian at all times.

%

As shown in Figure \ref{fig:feedback_architecture}, the objective is to design an adaptive feedback control law that drives the density matrix $\rho$ to a desired target state $\rho_\rmd$ without requiring prior knowledge of $H_0$, $H_1$, or $L_i$ in the control law. 
We assume that an estimate of the density matrix $\rho$ is available for feedback. 
In practice, the quantum state is not directly measurable in real time; rather, it must be reconstructed from measurement records using quantum filtering or observer-based techniques, such as those described in \cite{qamar2019observer}. 

\begin{remark}
    \textbf{Measurement and Observer Assumptions.}
\textit{The present formulation assumes access to an estimate of density $\rho$ for feedback and neglects the effect of measurement on the system dynamics. 
This corresponds to an ensemble-averaged or mean-field description of the open quantum system and serves as a first step toward evaluating model-free adaptive control at the level of LGKS dynamics. 
In general, a strict separation principle between quantum state estimation and feedback control does not hold for stochastic quantum systems with measurement backaction. 
Incorporating stochastic master equations and measurement-induced backaction into the feedback loop is an important direction for future work.} 
\end{remark}

\begin{figure}[h]
    \centering
    \resizebox{1\columnwidth}{!}{%
    \tikzset{every picture/.style={line width=1pt}} 
    
    \begin{tikzpicture}[x=0.75pt,y=0.75pt,yscale=-1,xscale=1, text width=2cm, align=center,minimum width=62pt, minimum height=32pt]
    
    \node[draw, fill=blue!20, rounded corners] (Box1) at (187+35,160) {$ F(\rho,\rho_{\rm{d}})$};
    \node[draw, fill=blue!20, rounded corners] (Box2) at (316+35,160) {Adaptive Controller};
    \node[draw, fill=blue!20, rounded corners] (Box3) at (435+35,160) {$\mathcal{L}\mathcal{G}\mathcal{K}\mathcal{S}$};
    
    \node[draw, fill=red!20, rounded corners] (Box4) at (420,256) {Homodyne Detection};
    \node[draw, fill=red!20, rounded corners] (Box5) at (300,256) {State Observer};

    \path (Box2.south) ++(-20,0) coordinate (Target1);
    \path (Box2.north) ++(30,0) coordinate (Target2);
    \path (Box1.south) ++(0,234-160) coordinate (Target3);

    \draw[-latex] (Box1.west)++(-40,0) -- (Box1.west) node [midway, above] {$\rho_{\rm_{d}}$};
    \draw (Box1.east)++(10,0) -- ++(0,40) -- ++(30,0) --(Target1);
    \draw[-latex] (Target2) -- ++(30,-30);
    \draw[-latex] (Box1.east) -- (Box2.west) node [midway, above] {$e$};
    \draw[-latex] (Box2.east) -- (Box3.west) node [midway, above] {$u$};
    \draw[-latex] (Box3.east) -- ++(40,0) node [midway, above] {$\rho$};
    \draw[-latex] (Box3.east)++(15,0) |- (Box4.east);
    \draw[-latex] (Box4.west) -- (Box5.east);
    \draw(Box5.west) -- (Target3);
    \draw[-latex] (Box1.south)++(0,74) -- (Box1.south) node [midway, xshift=-7] {$\rho$};
    
    
    \end{tikzpicture}
    }
    \caption{Adaptive feedback control to track desired density matrix $\rho_\rmd.$
    The estimator, shown in red, consisting of a homodyne detection and a state observer, provides the density matrix $\rho$ for feedback.}
    \label{fig:feedback_architecture}
\end{figure}

\subsection{Density Error}
As shown in \cite{jones2024controlling}, the \textit{Uhlmann–Jozsa fidelity} between two density matrices $\rho_1, \rho_2$ is defined as 
\begin{align}
    F(\rho_1,\rho_{2})
        &\isdef
            \left(\tr\left(\sqrt{\sqrt{\rho_1}
            \,\rho_2
            \,\sqrt{\rho_1}}\right)\right)^2
        .
        \label{eq:fidelity}
\end{align}
It follows from Proposition \ref{thm:uhlman_trace} that, if $\rho_1, \rho_2$ are the states of \eqref{eq:Lindblad_equation}, then $F(\rho_1,\rho_2) \in [0,1].$
Furthermore, $F(\rho_1,\rho_2) = 1$ if and only if $\rho_1 = \rho_2.$
Finally, as shown in \cite{jozsa1994fidelity}, in the case where $\rho_1, \rho_2 \in \mathbb{C}^{2 \times 2},$ the \textit{Uhlmann–Jozsa fidelity} is given by
\begin{align}
    F(\rho_1,\rho_2)
        &=
            \tr(\rho_1\rho_2) + 2\sqrt{\det (\rho_1)\det(\rho_2)}
\end{align}

In this paper, we define the scalar density-tracking error as
\begin{align}
    e
        &\isdef
            1 - F(\rho,\rho_{\rm{d}}) \in [0,1],
        \label{eq:error_fidelity}
\end{align}
where $\rho_{\rm{d}} \in \mathbb{C}^{2 \times 2}$ is the desired density.
This formulation of error ensures that the feedback control law is driven by a scalar variable, thereby avoiding any geometric inconsistencies in the simulation of the LGKS system. 
Note that the error $e = 0$ when $\rho$ is equal to $\rho_\rmd.$




\subsection{Adaptive PID Control}
In this work, we consider an adaptive proportional-integral-derivative (PID) controller, which can be written as
\begin{align}
    u(t)
        &=
            k_\rmp(t) e(t) + k_\rmi(t) \int_0^t e(\tau) \rmd\tau + k_\rmd(t) \dot e(t),
        \label{eq:control_structure}
\end{align}
where $k_\rmp(t), k_\rmi(t),$ and $k_\rmd(t)$ are the adaptive proportional, integral, and derivative gains. 
Note that the controller can be written in the regressor form as
\begin{align}
    u(t) = \Phi(t) \theta(t),
\end{align}
where 
\begin{align}
    \Phi(t)
        \isdef 
            \matl 
                e(t) &
                \gamma(t) & 
                \dot e(t)
            \matr, 
    \quad
    \theta(t)
        \isdef 
            \matl 
                k_\rmp(t) \\
                k_\rmi(t) \\
                k_\rmd(t)
            \matr,
\end{align}
and the integral state $\gamma(t)$ satisfies $\dot \gamma(t) = e(t).$

The objective of the adaptive controller is to update the gains in $\theta(t)$ online so that the density error $e(t)$ is driven toward zero. 
The gain update law is obtained using the continuous-time retrospective cost adaptive control (RCAC) algorithm described in Section \ref{sec:CTRCAC}. 
The PID structure is adopted due to its low-order structure, interpretability, and compatibility with linear parameterizations required by RCAC. 
More general controller parameterizations, such as full-state feedback, can be incorporated within the same adaptive framework when additional measurements or state features are available.



%
\section{Continuous Time Retrospective Cost Adaptive Control}
\label{sec:CTRCAC}

Consider a dynamic system 
\begin{align}
    \dot x(t) &= f(x(t),u(t)), \\
    y(t) &= g(x(t)),
\end{align}
where 
$x(t) $ is the state, 
$u(t) $ is the input, 
$y(t) $ is the measured output, and 
the {\rm vec}tor functions $f$ and $g$ are the dynamics and the output maps. 
Define the \textit{performance variable}
\begin{align}
    z(t) \isdef y(t) - r(t), 
    \label{eq:z_def}
\end{align}
where $r(t)$ is the exogenous reference signal. 
The objective is to design an adaptive output feedback control law to ensure that $z(t) \to 0.$

\begin{remark}
    \textbf{Specialization to LGKS fidelity tracking.}
\textit{In the present application, the system state $x(t)$ corresponds to the density matrix $\rho(t)$ governed by the LGKS dynamics \eqref{eq:Lindblad_equation} with controlled Hamiltonian \eqref{eq:H_controlled}. 
The measured output is chosen as a scalar function of the density matrix, namely the fidelity-based signal $y(t) \isdef F(\rho(t),\rho_\rmd)$, and the reference is $r(t) \isdef 1$. 
Consequently, the performance variable \eqref{eq:z_def} becomes
$z(t) = y(t) - r(t) = F(\rho(t),\rho_\rmd) - 1 = -e(t)$, 
where $e(t)$ is the density error defined in \eqref{eq:error_fidelity}. 
This choice yields a scalar performance variable that is compatible with the linearly parameterized controller \eqref{eq:u_para}.}
\end{remark}

Consider a linearly parameterized control law 
\begin{align}
    u(t) 
        =
            \Phi(t) \theta(t), 
    \label{eq:u_para}
\end{align}
where the regressor matrix $\Phi(t) $ contains the measured data and 
the vector $\theta(t) $  contains the controller gains to be optimized. 
Various linear parameterizations of MIMO controllers are described in \cite{goel_2020_sparse_para}.

Next, using \eqref{eq:z_def}, define the \textit{retrospective performance}
\begin{align}
    \hat z(t) 
        &\isdef
            z(t) + \Phi_\rmf(t) \hat \theta(t) -u_\rmf(t),
\end{align}
where $\hat \theta(t)$ is the controller gain to be optimized and  
the filtered regressor $\Phi_\rmf(t)$ and the filtered control $u_\rmf(t)$ are defined as
\begin{align}
    \Phi_\rmf(t) &\isdef G_\rmf(s) \left[ \Phi(t) \right] ,
    \\
    u_\rmf(t) &\isdef G_\rmf(s) \left[ u(t) \right] ,
\end{align}
where $G_\rmf(s)$ is a dynamic filter.


The filtered signals $\Phi_\rmf(t)$ and $u_\rmf(t)$ are computed in the time domain, as shown below. 
Let $(A_\rmf,B_\rmf,C_\rmf,D_\rmf)$ be a realization of $G_\rmf(s).$
Then, 
\begin{align}
    \dot x_\Phi &= A_\rmf x_\phi + B_\rmf \Phi, \\
    \Phi_\rmf &= C_\rmf x_\phi + D_\rmf \Phi,
\end{align}
and
\begin{align}
    \dot x_u &= A_\rmf x_u + B_\rmf u, \\
    u_\rmf &= C_\rmf x_u + D_\rmf u.
\end{align}

\begin{remark}
    \textbf{Filter states and notation.}\textit{
The states $x_\phi$ and $x_u$ appearing in the filter realizations below are internal states of the dynamic filter $G_\rmf(s)$ driven by $\Phi(t)$ and $u(t)$, respectively. 
These filter states are auxiliary variables used only to compute $\Phi_\rmf(t)$ and $u_\rmf(t)$ and are not related to the LGKS state $\rho(t)$.}
\end{remark}

Next, define the retrospective cost 
\begin{align}
    J(t, \hat \theta)
        &=
            \int_{0}^{t} 
            \Big(
                e^{-\lambda(t-\tau)}
                \hat z(\tau)^\rmT R_z \hat z(\tau) 
                \nn \\ &\quad
                +
                (\Phi(\tau)   \hat \theta)^\rmT R_{u} (\Phi(\tau)  \hat \theta)
            \Big)
            \rmd \tau
            +
            e^{-\lambda t}
            \hat \theta^T R_\theta \hat \theta,
    \label{eq:J_RCAC}
\end{align}
where $R_z,
R_u,$ and $R_\theta$ are positive definite weighting matrices of appropriate dimensions, and 
$\lambda > 0$ is the exponential forgetting factor. 

The retrospective cost \eqref{eq:J_RCAC} defines an online performance objective for updating the parameter vector $\theta(t)$ based on measured data. 
We emphasize that minimizing \eqref{eq:J_RCAC} does not, in general, imply closed-loop stability or guarantee that $z(t)\to 0$ for nonlinear systems. 
In particular, for the LGKS dynamics considered in this paper, the tracking behavior reported in Section \ref{sec:numerical_simulation} is empirical and depends on the system, the controller parameterization, and the choice of hyperparameters.

\begin{proposition}
    Consider the cost function $J(t, \hat \theta)$ given by \eqref{eq:J_RCAC}.
    For all $t \ge 0, $ define the minimizer of $J(t, \hat \theta)$ by
    \begin{align}
        \theta(t)
            \isdef 
                \underset{\hat \theta \in \BBR^{l_\theta}}{\operatorname{argmin}} 
                J(t, \hat \theta).
    \end{align}
    Then, for all $t \ge 0, $ the minimizer satisfies
    \begin{align}
        \dot \theta
        &=
            -P \Phi_\rmf^\rmT R_z (z+ \Phi_\rmf  \theta -u_\rmf  )
            -P \Phi^\rmT R_{u} \Phi \theta, 
        \label{eq:theta_eq}
        \\
        \dot P
        &=
            \lambda P
            -P
            \left(
            \Phi_\rmf ^\rmT 
                    R_z
                    \Phi_\rmf  
                    +
                    \Phi^\rmT R_{u} \Phi
            \right)
            P
    \end{align}
    where $P(0) = R_\theta^{-1}$ and $\theta(0) = 0.$
    
\end{proposition}

\begin{proof}
    Note that the cost function \eqref{eq:J_RCAC} can be written as
\begin{align}
    J(t, \hat \theta) 
        &=
            \hat \theta^T A(t) \hat \theta + 2 \hat \theta^T b(t) + c(t),
    \label{eq:J_reformulated}
\end{align}
where
\begin{align}
    A(t) 
        &\isdef
            \int_0^{t} 
                 e^{-\lambda(t-\tau)}
            \left(
                    \Phi_\rmf ^\rmT 
                    R_z
                    \Phi_\rmf  
                    +
                    \Phi^\rmT R_{u} \Phi 
            \right)
            \rmd \tau 
            +
             e^{-\lambda t} R_\theta,
    \nn
    \\
    b(t) 
        &\isdef
            \int_0^{t}
             e^{-\lambda(t-\tau)}
            \Phi_\rmf^\rmT
                R_z 
                (z-u_\rmf) \rmd \tau, 
    \nn \\
    c(t) 
        &\isdef
            \int_0^{t}
             e^{-\lambda(t-\tau)}
            (z-u_\rmf)^\rmT
                R_z 
                (z-u_\rmf)  
                \rmd \tau.
    \nn
\end{align}
The matrices $A(t)$ and $b(t)$ satisfy
\begin{align}
    \dot A
        &=
            -\lambda A +
            \Phi_\rmf ^\rmT 
                    R_z
                    \Phi_\rmf  
                    +
                    \Phi^\rmT R_{u} \Phi,  \nn 
    \\
    \dot b
        &=
            -\lambda b +
            \Phi_\rmf^\rmT
                R_z 
                (z-u_\rmf) , \nn 
\end{align}
where $A(0) = R_\theta$ and $b(0) = 0.$


Next, define, for all $t \ge 0,$ $P(t) \isdef A(t)^{-1}.$
Using the fact that $P(t) A(t) = I,$ it follows that
\begin{align}
    \dot P(t) = - P(t) \dot A(t) P(t), \nn
\end{align}
and thus 
\begin{align}
    \dot P 
        =
            \lambda P
            -P
            \left(
            \Phi_\rmf ^\rmT 
                    R_z
                    \Phi_\rmf  
                    +
                    \Phi^\rmT R_{u} \Phi
            \right)
            P. \nn
\end{align}

Finally, note that the minimizer of \eqref{eq:J_reformulated} is given by
\begin{align}
    \theta(t)
        =
            -A(t)^{-1} b(t)
        =
            -P(t) b(t), \nn
\end{align}
and thus 
\begin{align}
    \dot \theta(t) 
        &=
            -\dot P b - P \dot b
        \nn 
        \\
        &=
            -P
            \left(
            \Phi_\rmf ^\rmT 
                    R_z
                    \Phi_\rmf  
                    +
                    \Phi^\rmT R_{u} \Phi
            \right) P
            P\inv \theta
            \nn \\ &\quad -
            P \Phi_\rmf^\rmT R_z (z-u_\rmf)
        \nn \\
        &=
            -P \Phi_\rmf^\rmT R_z (z+ \Phi_\rmf  \theta -u_\rmf  )
            -P \Phi^\rmT R_{u} \Phi \theta.
        \nn
\end{align}
\end{proof}


The control law is thus
\begin{align}
    u(t) = \Phi(t) \theta(t), 
\end{align}
where $\theta(t)$ is given by \eqref{eq:theta_eq}.

For the adaptive PID structure in Section \ref{sec:problem_formulation}, the regressor $\Phi(t)$ and parameter vector $\theta(t)$ in \eqref{eq:u_para} correspond to the fidelity-based features and gains defined in \eqref{eq:control_structure}. 
Thus, the update laws \eqref{eq:theta_eq} provide an online method to adjust $(k_\rmp(t),k_\rmi(t),k_\rmd(t))$ using only the measured fidelity-based error signal and its integral/derivative features, without requiring explicit knowledge of $H_0$, $H_1$, or $L_i$ in the control law.

%


\section{Numerical Simulation}
\label{sec:numerical_simulation}

In this paper, we consider a two-level LGKS system \eqref{eq:Lindblad_equation} with a single dissipative channel, that is, $m = 1$. 
The control system architecture is shown in Fig.~\ref{fig:feedback_architecture}. 
Note that, for conceptual clarity, the observer is included in the block diagram; however, in the numerical studies, we assume access to an estimate of the density matrix $\rho$ for feedback computation. 
In practice, the quantum state is not directly measurable in real time and must be reconstructed from measurement records using quantum filtering or observer-based techniques; a representative observer design is described in \cite{qamar2019observer}. 
In this preliminary study, the observer dynamics and measurement backaction are not explicitly modeled, and the focus is on evaluating the closed-loop behavior of the adaptive controller under idealized assumptions about state access.

To simulate the LGKS equation, we set 
\begin{align}
    H_0
        =
            0.5
                \matl
                    1 & 0\\
                    0 & -1
                \matr,
    H_1
        =
            0.5
                \matl
                    0 & 1\\
                    1 & 0
                \matr,
    L_1
        =
            \matl
                0 & 1\\
                0 & 0
            \matr,
\end{align}
where $H_0$ and $H_1$ are \textit{Pauli matrices} that, as described in \cite{d2001optimal,kamakari2022digital}, satisfy the controllability conditions in the problem of controlling a two-level quantum system, 
and $L_1$ is selected as a non-Hermitian matrix, as suggested in \cite{verstraete2009quantum,schlimgen2022quantum}.

The initial density matrix is chosen as
\begin{align}
    \rho_0
        &=
            \matl
                0.4 & 0.1+0.3\imath\\
                0.1-0.3\imath & 0.6
            \matr,
    \label{eq:rho_0}
\end{align}
which is a randomly generated Hermitian, positive semidefinite matrix with unit trace. 
Following \cite{moses2015creation, xu2021multicriticality}, we consider both low-entropy and high-entropy target states $\rho_\rmd$, where the \textit{von Neumann entropy} of a quantum state $\rho$ is defined as
\begin{align}
    S(\rho) \isdef -\tr(\rho\,\rm{ln}(\rho)).
\end{align}


Achieving or maintaining low-entropy states is often a central objective in quantum information processing, since low entropy corresponds to high purity and coherence, which are essential resources for quantum computation, communication, and sensing. 
Because environmental interactions and decoherence naturally increase entropy, feedback control is required to counteract this drift and preserve useful quantum states. 
Conversely, the ability to track high-entropy target states is also relevant in many physical scenarios, since open quantum systems are often driven toward mixed or thermal equilibrium states by their environment.


\subsection{Hyperparameter Tuning}

First, we specify the low-entropy target state
\begin{align}
    \rho_\rmd
        =
            \matl
                0.8571  & 0.2857 + 0.1429\imath\\
                 0.2857 - 0.1429\imath &  0.1429
            \matr.
    \label{eq:rho_d_LE}
\end{align}
Note that $S(\rho_\rmd)  = 0.1013,$ and thus, this desired state is a low-entropy state. 
Furthermore, $\rho_\rmd,$ given by \eqref{eq:rho_d_LE}, is an equilibrium point of \eqref{eq:Lindblad_equation} in the case where $u(t) \equiv 1.$

With the desired density $\rho_\rmd$ given by \eqref{eq:rho_d_LE} as the reference for the controller, we tune the RCAC hyperparameters. 
In particular, we set $R_z = 1, $ $R_u = 1,$ and $\lambda = 0.01,$ and optimize $P_0$ and the filter coefficient $\beta,$ where $G_\rmf(s) = \dfrac{1}{s+\beta}.$ 
We restrict the optimization to $P_0$ and $\beta,$ as numerical experiments consistently indicate that these parameters are the most dominant factors in the closed-loop performance.

\begin{figure}[!ht]
    \centering
    \includegraphics[width=\linewidth]{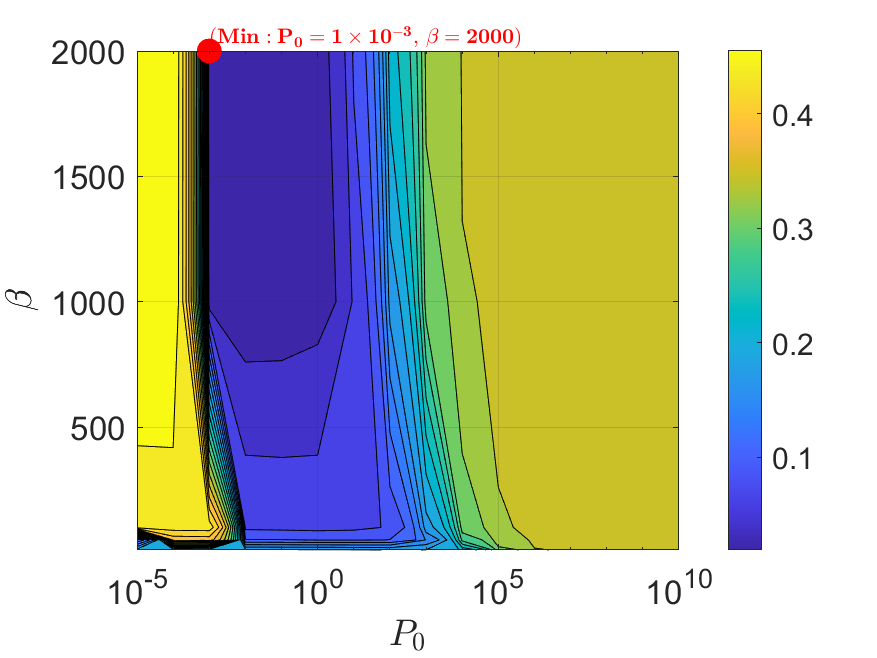}
    \caption{Contour plot of $J_h$ for various values of $(P_0, \beta).$
    The pair $(P_0, \beta)$ that minimizes $J_h$ is used in all numerical simulations.}
    \label{fig:tunning_figure_scalar_case}
\end{figure}

For each pair of $(P_0, \beta) \in \{ (10^{-5}, 10^{-4}, \ldots, 10^{10})  \times (0, 1, 2, 5, 100, 2000)\},$ the RCAC algorithm is applied to update the adaptive PID controller to drive the quantum state from $\rho_0$, given by \eqref{eq:rho_0}, to $\rho_\rmd$, given by \eqref{eq:rho_d_LE}. 
Each simulation is run for 200 seconds. 
At the end of each simulation, the hyperparameter cost 
\begin{align}
    J_h
        \isdef
            \int_{190}^{200} |z(\tau)| \rmd \tau
\end{align}
is computed. 
Figure \ref{fig:tunning_figure_scalar_case} shows the contour plot of $J_h$ for various values of $(P_0, \beta).$ 
The tuning pair with the minimum value of $J_h$ is selected for all numerical simulations presented in the remainder of the paper. 
In particular, in all the examples, we set $P_0 = 10^{-3}\times I_3$ and $\beta = 2000.$ 
We emphasize that the hyperparameter tuning procedure outlined above is coarse and heuristic in nature, and is intended to identify a representative near-optimal parameter set for the examples considered here.

\subsection{Low-Entropy Density Tracking}

\begin{figure}[ht]
    \centering
    \includegraphics[width=0.9\columnwidth]{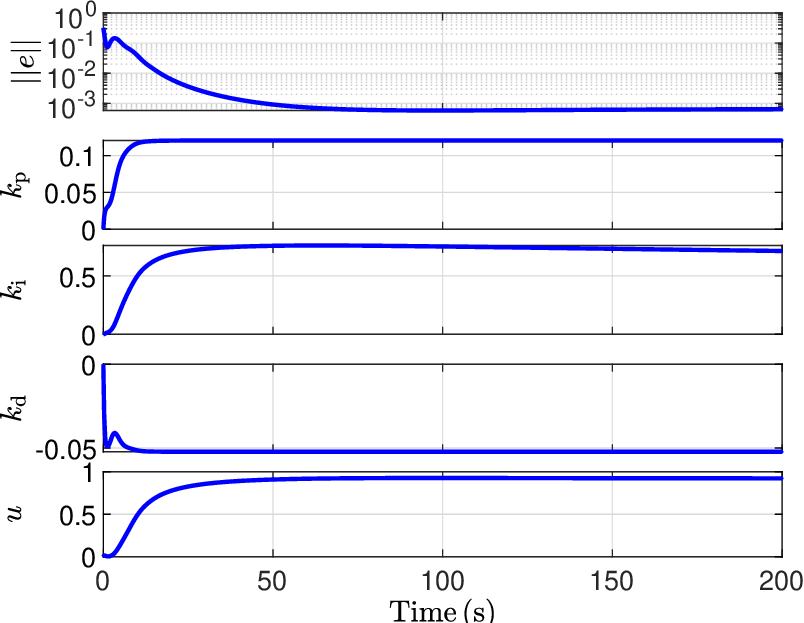}
    \caption{
    \textbf{Low-entropy density tracking.}
    The first subplot shows the absolute value of the density error $e(t)$ between the system state and the desired state on a logarithmic scale, 
    the next three subplots show the adaptive gains $k_\rmp,$ $k_\rmi,$ and $k_\rmd,$ updated by RCAC, and the last subplot shows the control $u$ generated by the adaptive controller.}
    \label{fig:PID_RCAC_forgetting_factor_Quantum_system_scalar_u_parameters_and_error_paper}
\end{figure}

\begin{figure}[ht]
    \centering
    \includegraphics[width=0.9\columnwidth]{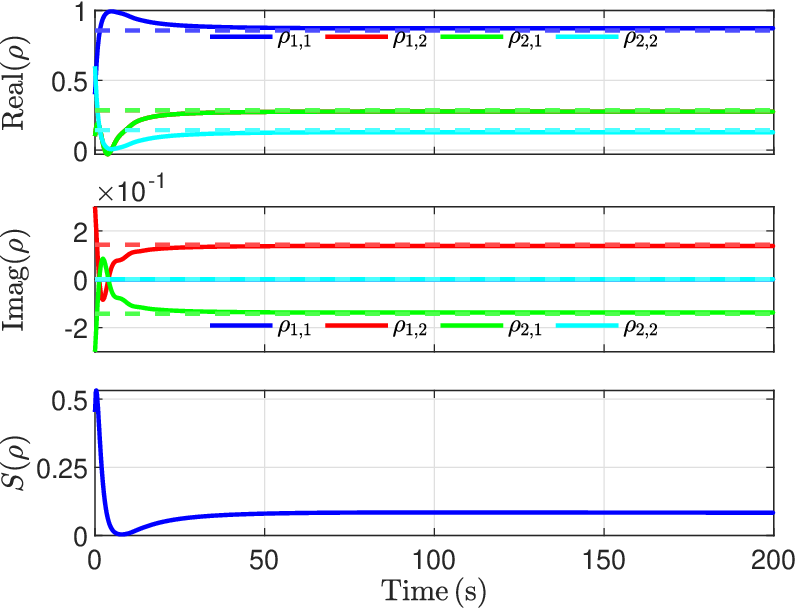}
    \caption{
    \textbf{Low-entropy density tracking.}
    The first two subplots show the real and imaginary parts of the quantum state $\rho(t).$ 
    The desired state components are shown with dashed lines, while the system response is shown with solid lines of the corresponding color. 
    The third subplot shows the entropy of the quantum state.}
    \label{fig:PID_RCAC_forgetting_factor_Quantum_system_scalar_u_rho_and_rho_desired_forget}
\end{figure}

\begin{figure}[ht]
    \centering
    \includegraphics[width=0.9\columnwidth]{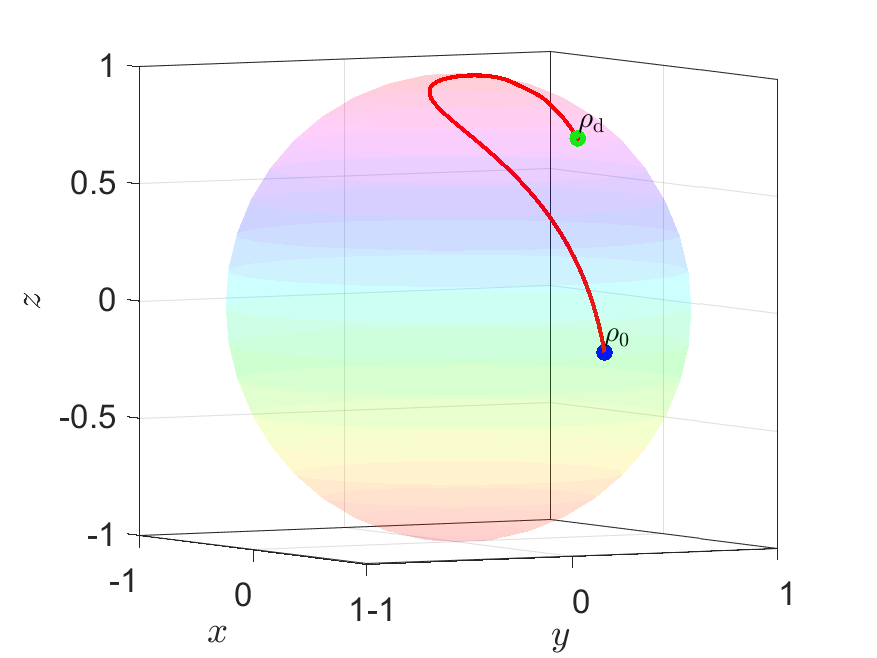}
    \caption{
    \textbf{Low-entropy density tracking.}
    Trajectory of the quantum state of the LGKS system \eqref{eq:Lindblad_equation} on the Bloch sphere. 
    Note that the trajectory remains inside the Bloch ball, consistent with mixed-state evolution.}
    \label{fig:PID_RCAC_forgetting_factor_Quantum_system_scalar_bloch_sphere}
\end{figure}

With the selected RCAC hyperparameters, 
Fig.~\ref{fig:PID_RCAC_forgetting_factor_Quantum_system_scalar_u_parameters_and_error_paper} shows the absolute value of the density error $e(t)$ between the detected system state and the desired low-entropy state on a logarithmic scale, the adaptive gains $k_\rmp,$ $k_\rmi,$ and $k_\rmd,$ updated by RCAC, and the control $u$ generated by the adaptive controller. 
Fig.~\ref{fig:PID_RCAC_forgetting_factor_Quantum_system_scalar_u_rho_and_rho_desired_forget} shows the real and imaginary parts of the quantum state $\rho(t)$ and the von Neumann entropy of the quantum state. 
Fig.~\ref{fig:PID_RCAC_forgetting_factor_Quantum_system_scalar_bloch_sphere} shows the trajectory of the quantum state of the LGKS system \eqref{eq:Lindblad_equation} on the Bloch sphere. 
Note that the quantum state $\rho$ satisfies $x = 2\,\rm{Real}(\rho_{1,2}),$ $y = 2\,\rm{Imag}(\rho_{1,2}),$ and $z = \rho_{1,1} - \rho_{2,2}$ on the Bloch sphere \cite{wie2020two}.

\subsection{High-Entropy Density Tracking}

Next, we consider the problem of driving the quantum state $\rho(t)$ to the high-entropy target
\begin{align}
    \rho_\rmd
        =
            \matl
                0.5168  & 0.0971 + 0.0460i\\
                0.0971 - 0.0460i &  0.4832
            \matr.
    \label{eq:rho_d_HE}
\end{align}
Note that $S(\rho_\rmd)  = 0.6693,$ and thus, this desired state is a high-entropy state. 
Furthermore, $\rho_\rmd,$ given by \eqref{eq:rho_d_HE}, is an equilibrium point of \eqref{eq:Lindblad_equation} in the case where $u(t) \equiv 10.$

\begin{figure}[ht]
    \centering
    \includegraphics[width=0.9\columnwidth]{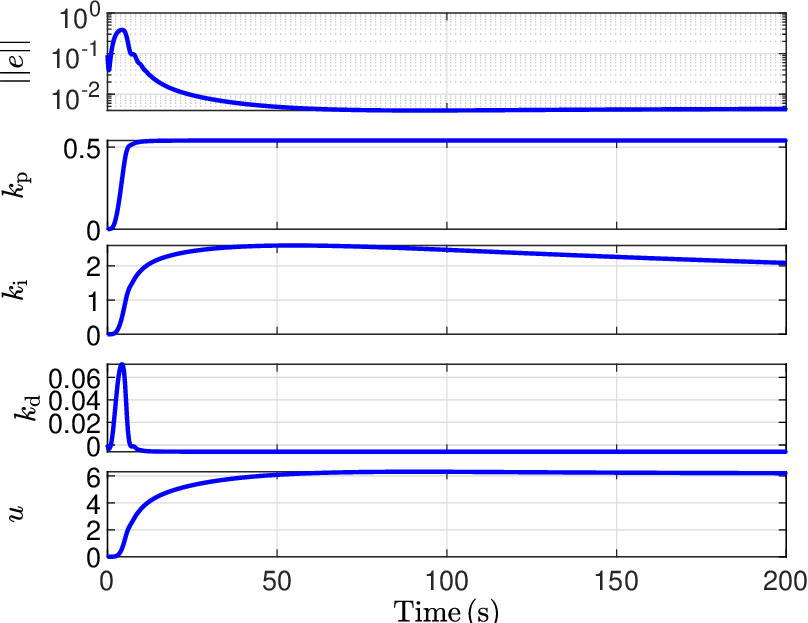}
    \caption{
    \textbf{High-entropy density tracking.}
    The first subplot shows the absolute value of the density error $e(t)$ between the system state and the desired state on a logarithmic scale,
    the next three subplots show the adaptive gains $k_\rmp,$ $k_\rmi,$ and $k_\rmd,$ updated by RCAC, and the last subplot shows the control $u$ generated by the adaptive controller.}
    \label{fig:PID_RCAC_Quantum_Parameters_and_errors_high_entropy}
\end{figure}

\begin{figure}[!ht]
    \centering
    \includegraphics[width=0.9\columnwidth]{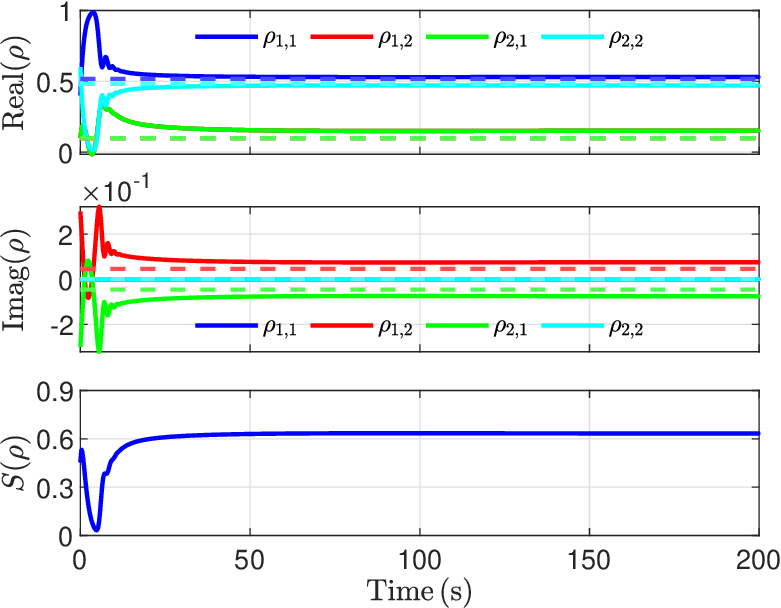}
    \caption{
    \textbf{High-entropy density tracking.}
    The first two subplots show the real and imaginary parts of the quantum state $\rho(t).$ 
    The desired state components are shown with dashed lines, while the system response is shown with solid lines of the corresponding color. 
    The third subplot shows the entropy of the quantum state.}
    \label{fig:PID_RCAC_Quantum_Parameters_rho_and_rho_desired_high_entropy}
\end{figure}

\begin{figure}[!ht]
    \centering
    \includegraphics[width=0.9\columnwidth]{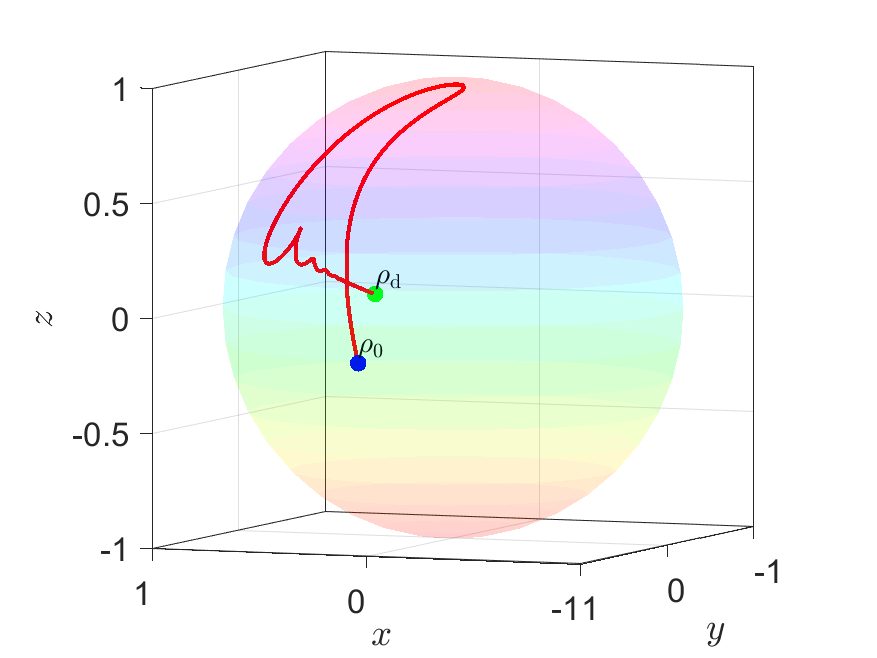}
    \caption{
    \textbf{High-entropy density tracking.}
    Trajectory of the quantum state of the LGKS system \eqref{eq:Lindblad_equation} on the Bloch sphere. 
    Note that the trajectory remains inside the Bloch ball, consistent with mixed-state evolution.}
    \label{fig:PID_RCAC_Quantum_bloch_high_entropy}
\end{figure}

Figure~\ref{fig:PID_RCAC_Quantum_Parameters_and_errors_high_entropy} shows the absolute value of the density error $e(t)$ between the detected system state and the desired high-entropy state on a logarithmic scale, the adaptive gains $k_\rmp,$ $k_\rmi,$ and $k_\rmd,$ updated by RCAC, and the control input $u$ generated by the adaptive controller. 
Figure~\ref{fig:PID_RCAC_Quantum_Parameters_rho_and_rho_desired_high_entropy} shows the real and imaginary parts of the quantum state $\rho(t)$ and the von Neumann entropy of the quantum state. 
Figure~\ref{fig:PID_RCAC_Quantum_bloch_high_entropy} shows the trajectory of the quantum state of the LGKS system \eqref{eq:Lindblad_equation} on the Bloch sphere.

\subsection{Shot Noise}

Finally, we revisit the low-entropy target $\rho_\rmd$ given by \eqref{eq:rho_d_LE} and examine the effect of measurement noise on the feedback loop. 
Shot noise in the feedback detection is modeled as additive white noise, following \cite{niebauer1991nonstationary}. 
Specifically, at each discrete sampling instant, a perturbation $N \in \mathbb{C}^{2 \times 2}$ is drawn with entries distributed as $\mathcal{N}(0,\sigma^{2})$, and the noisy state is constructed as
\begin{equation}
    \rho_{\text{noisy}} = \frac{\rho + \tfrac{1}{2}(N + N^{\dagger})}{\tr\!\left(\rho + \tfrac{1}{2}(N + N^{\dagger})\right)}.
\end{equation}
Hermiticity is enforced by symmetrization, and physical validity is preserved by trace normalization, ensuring that $\rho_{\text{noisy}}$ remains a valid density matrix under white noise of strength $\sigma = 5 \times 10^{-3}$. 

\begin{remark}
    \textbf{Noise model and feedback-only corruption.}
\textit{The noisy state $\rho_{\text{noisy}}$ is used only in the feedback path for computing the fidelity-based error signal. 
The underlying system evolution is still governed by the deterministic LGKS dynamics driven by the true state $\rho(t)$. 
Thus, the injected noise models the corruption of the measurement or state estimate rather than physical noise acting directly on the quantum dynamics.}
\end{remark}

\begin{figure}[ht]
    \centering
    \includegraphics[width=0.9\columnwidth]{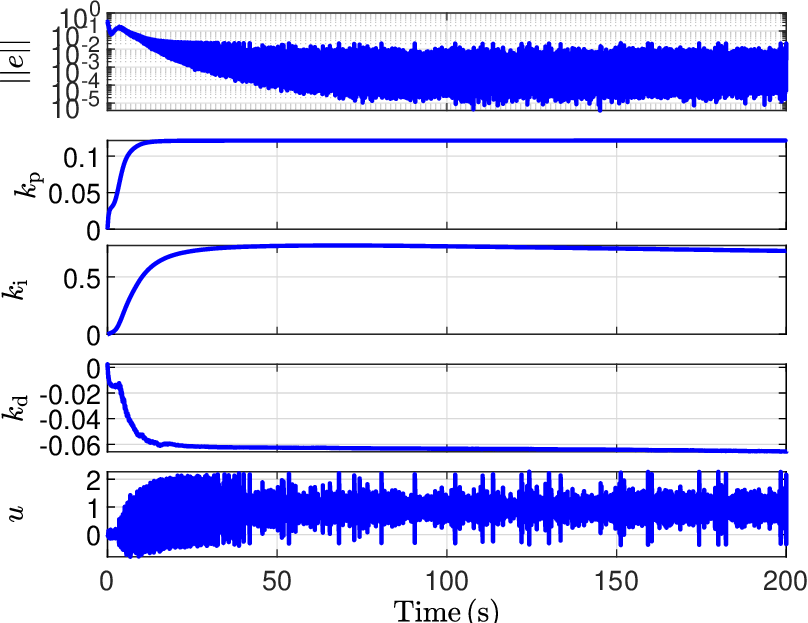}
    \caption{
    \textbf{Shot noise case tracking.}
    The first subplot shows the absolute value of the density error $e(t)$ between the detected system state with shot noise and the desired state on a logarithmic scale,
    the next three subplots show the adaptive gains $k_\rmp,$ $k_\rmi,$ and $k_\rmd,$ updated by RCAC, and the last subplot shows the control $u$ generated by the adaptive controller.}
    \label{fig:PID_RCAC_forgetting_factor_Quantum_system_scalar_u_parameters_and_error_paper_noise}
\end{figure}

\begin{figure}[ht]
    \centering
    \includegraphics[width=0.9\columnwidth]{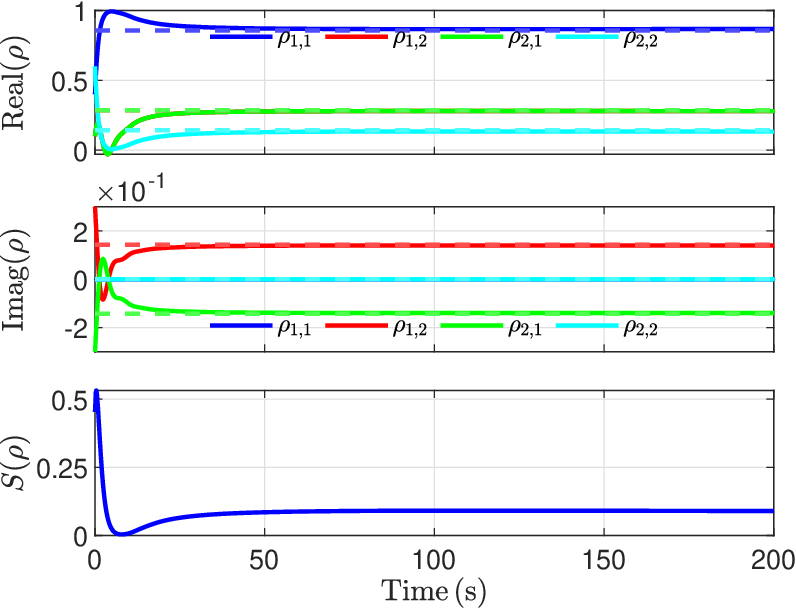}
    \caption{
    \textbf{Shot noise case tracking.}
    The first two subplots show the real and imaginary parts of the quantum state $\rho(t).$ 
    The desired state components are shown with dashed lines, while the system response is shown with solid lines of the corresponding color. 
    The third subplot shows the entropy of the quantum state.}
    \label{fig:PID_RCAC_forgetting_factor_Quantum_system_scalar_u_rho_and_rho_desired_forget_noise}
\end{figure}

\begin{figure}[ht]
    \centering
    \includegraphics[width=0.9\columnwidth]{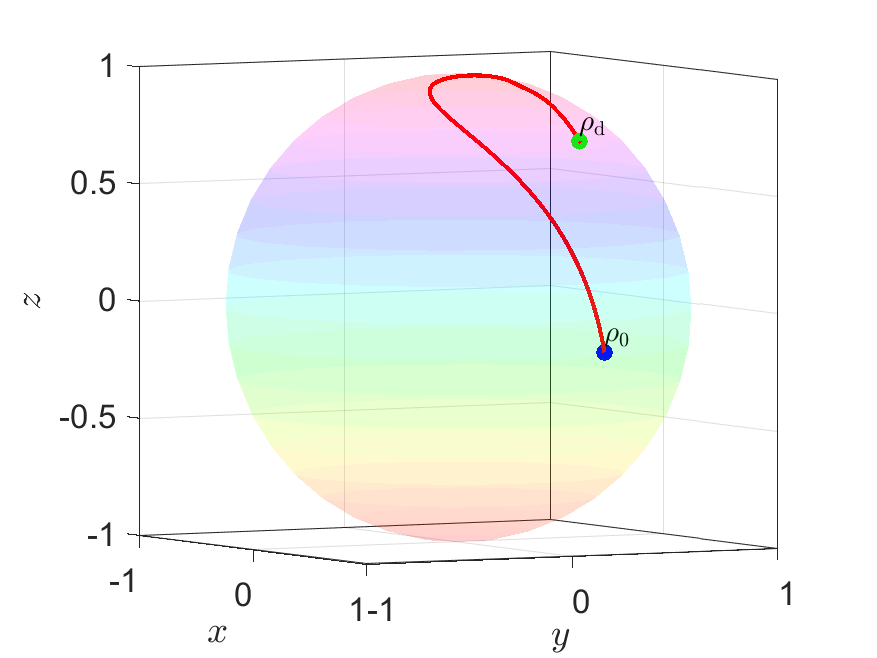}
    \caption{
    \textbf{Shot noise case tracking.}
    Trajectory of the quantum state of the LGKS system \eqref{eq:Lindblad_equation} on the Bloch sphere. 
    Note that the trajectory remains inside the Bloch ball, consistent with mixed-state evolution.}
    \label{fig:PID_RCAC_forgetting_factor_Quantum_system_scalar_bloch_sphere_noise}
\end{figure}

Figure~\ref{fig:PID_RCAC_forgetting_factor_Quantum_system_scalar_u_parameters_and_error_paper_noise} shows the absolute value of the density error $e(t)$ between the detected system state with shot noise and the desired state on a logarithmic scale, the adaptive gains $k_\rmp,$ $k_\rmi,$ and $k_\rmd,$ updated by RCAC, and the control $u$ generated by the adaptive controller. 
Figure~\ref{fig:PID_RCAC_forgetting_factor_Quantum_system_scalar_u_rho_and_rho_desired_forget_noise} shows the real and imaginary parts of the quantum state $\rho(t)$ and the von Neumann entropy of the quantum state. 
Figure~\ref{fig:PID_RCAC_forgetting_factor_Quantum_system_scalar_bloch_sphere_noise} shows the trajectory of the quantum state of the LGKS system \eqref{eq:Lindblad_equation} on the Bloch sphere. 

The three examples described above suggest that the adaptive PID controller updated by the RCAC algorithm can effectively track a desired density matrix $\rho_\rmd$ in the presence of moderate measurement noise, without requiring a model of the underlying LGKS dynamics in the control law.

\section{Conclusions}
\label{sec:conclusions}

This paper presented an online learning-based adaptive control framework for tracking a desired density matrix in a two-level Lindblad--Gorini--Kossakowski--Sudarshan quantum system, in which the feedback control law does not require prior knowledge of the system Hamiltonian or dissipative operators. 
In particular, the adaptive controller is based on the retrospective cost adaptive control (RCAC) algorithm, and a continuous-time formulation is employed to align with the LGKS dynamics. 
An adaptive PID controller structure is used to generate the control signal while preserving the geometric structure of the LGKS evolution. 
Uhlmann--Jozsa fidelity is adopted as a scalar performance metric to drive the adaptation. 
Numerical simulations demonstrate effective tracking of both low-entropy and high-entropy target states under idealized state-access assumptions and in the presence of moderate measurement noise.

Future work will extend this framework to matrix-valued control inputs and alternative quantum error metrics, such as coherence-based measures, enabling richer controller parameterizations and evaluation on higher-dimensional quantum systems. 
Additional directions include incorporating stochastic measurement backaction, discrete-time feedback implementations, and theoretical analysis of closed-loop performance guarantees.

\renewcommand{\thesection}{\Alph{section}}
\section*{Appendix}
\setcounter{section}{0}

\section{Some Useful Facts}
\label{sec:preliminaries}

This appendix reviews some useful definitions and facts about Hermitian matrices used in the LKGS system. 


\begin{definition}
    Let $A \in \BBC^{n \times n}.$
    $A$ is a Hermitian matrix if $A = \overline{A^\rmT},$
    where $\overline{A^\rmT}$ is the conjugate transpose of $A.$
    
\end{definition}

The conjguate transpose of $A$ is denoted by $A^\rmH.$
If $A$ is Hermitian, then $A = A^\rmH.$

\begin{definition}
    Let $A, B \in \BBC^{n \times n}.$
    Then, the commutator $[A,B]$ is defined as 
    \begin{align}
        [A,B] \isdef AB - BA,
    \end{align}
    and the anticommutator $\{A,B \}$ is defined as 
    \begin{align}
        \{A,B\} \isdef AB + BA.
    \end{align}
\end{definition}


\begin{fact}
    \label{thm:herm_preserve_1}
    Let $A, B \in \BBC^{n \times n}$ be Hermitian.
    Then, $-\imath[A,B]$ is Hermitian.
\end{fact}
\begin{proof}
    Note that 
    $(-i[A,B])^\rmH  
    = (-\imath)^\rmH [A,B]^\rmH  
    = \imath\left(AB - BA\right)^\rmH  
    = \imath(A^\rmH B^\rmH  - B^\rmH A^\rmH ) 
    = \imath(BA - AB) 
    = -\imath(AB - BA) 
    = -\imath[A,B]$. 
    Thus, $-\imath[A,B]$ is Hermitian.   
\end{proof}

\begin{fact}
    \label{thm:herm_preserve_2}
    Let $U \in \BBC^{n \times n}$ be Hermitian. 
    Then, $U ^\rmH U$ and $U U^\rmH$ are Hermitian. 
\end{fact}
\begin{proof}
    Note that $(U^\rmH U)^\rmH  = U^\rmH (U^\rmH )^\rmH  = U^\rmH U,$ which implies that $UU^\rmH $ is Hermitian. 
    A similar argument proves that $U^\rmH U$ is also Hermitian. 
\end{proof}

\begin{fact}
    \label{thm:herm_preserve_3}
    Let $A, U \in \BBC^{n \times n}$ be Hermitian.
    Then, $UAU^\rmH $ is Hermitian.
\end{fact}
\begin{proof}
    Note that $(UAU^\rmH )^\rmH  = (U^\rmH )^\rmH (UA)^\rmH  = UA^\rmH U^\rmH  = UAU^\rmH ,$ which implies that $UAU^\rmH $ is Hermitian.
\end{proof}

\begin{fact}
    \label{thm:trace_peservation}
    Let $A, B, C \in \BBC^{n \times n}$ be Hermitian. 
    Suppose $B$ satisfies
    \begin{align}
        \dot{B} 
            =
                -\imath[A,B] + \sum_{i=1}^N\left(CB C^\rmH  - \dfrac{1}{2}\{C^\rmH C, B\}\right).
            \label{eq:theorem_trace_preservation}
    \end{align}
    Then, for all $t \geq 0,$ $\tr B $ is a constant.

\end{fact}
\begin{proof}
    %
    Note that $\dfrac{\rmd }{\rmd t} \tr(B)= \tr(\dot{B})$. 
    Next, note that $\tr(-\imath[A,B]) = -\imath\tr(AB - BA) = -\imath\tr(AB) + \imath\tr(BA) = -\imath\tr(AB) + \imath\tr(AB) = 0.$ 
    Using the cyclic property of the trace, it follows that 
    \begin{align*}
        &\tr\left(\sum_i^{N}\left(CBC^\rmH  - \dfrac{1}{2}(C^\rmH CB + BC^\rmH C)\right)\right) 
            \nn \\
            &=
                \sum_i^{N}\tr(CBC^\rmH ) -\dfrac{1}{2}\sum_i^N\tr(C^\rmH CB) - \dfrac{1}{2}\sum_i^N\tr(C^\rmH CB) 
            \nn \\
            &=
                \sum_i^{N}\tr(C^\rmH CB) - \sum_i^{N}\tr(C^\rmH CB) 
            =
                0,
    \end{align*}
    which completes the proof.     
\end{proof}

\begin{proposition}
    \label{thm:uhlman_trace}
    Let $\rho \in \mathbb{C}^{n \times n}$ and $\sigma \in \mathbb{C}^{n \times n}$ be two Hermitian PSD matrices with traces equal to $1$. Then,
    \begin{align}
    1
        \geq
            \left(\tr\left(\sqrt{\sqrt{\rho}
            \, \sigma 
            \,\sqrt{\rho}}\right)\right)^2 \in \mathbb{R}
                \geq 0.\nn   
    \end{align}
\end{proposition}
\begin{proof}
    Since $\rho$ and $\sigma$ are PSD Hermitian. Then, $\rho^{\rmH} = \rho$, $\sigma^{\rmH} = \sigma$, $x^{\rmH}\rho x\geq 0$ and $x^{\rmH}\sigma x\geq 0$ for all $x \in \mathbb{C}.$ Furthermore, the unique square root of a PSD Hermitian matrix is also Hermitian and PSD \cite[p~440]{horn2012matrix}. This implies that,
        $\left(\sqrt(\rho)\sigma\sqrt{\rho}\right)^{\rmH}
            =
                \sqrt{\rho}^{\rmH}\left(\sqrt(\rho)\sigma\right)^{\rmH}
                    =
                        \sqrt{\rho}^{\rmH}\sigma^{\rmH}\sqrt{\rho}^{\rmH}
                            =\sqrt{\rho}\sigma\sqrt{\rho}.$
    Define, $y \isdef \sqrt{\rho}x \in \mathbb{C}^n.$ Then, $x^{\rmH}\sqrt{\rho}\sigma\sqrt{\rho}x = y^{\rmH}\sigma y \geq 0,$ which implies that $\sqrt{\rho}\sigma\sqrt{\rho} \geq0.$
    Therefore, $\left(\sqrt(\rho)\sigma\sqrt{\rho}\right)$ is Hermitian, and $\sqrt{\left(\sqrt(\rho)\sigma\sqrt{\rho}\right)}$ is Hermitian and unique. From the Cauchy-Schwarz inequality and the cyclic property for traces, $|\tr(\sqrt{\sqrt{\rho}\sigma\sqrt{\rho}})|^2 = |\tr(\sqrt{\sqrt{\rho}\sqrt{\rho}\sigma})|^2 = |\tr(\sqrt{\rho\sigma})|^2 \leq \tr(\rho)\tr(\sigma) \leq 1$. Thus,
    $1
        \geq
            \left(\tr\left(\sqrt{\sqrt{\rho}
            \,\rho_{\rm{d}}
            \,\sqrt{\rho}}\right)\right)^2 \in \mathbb{R} 
                \geq 0.$
\end{proof}

\printbibliography

\end{document}